\documentclass[12pt]{iopart}
\usepackage[T1]{fontenc}
\usepackage[ansinew]{inputenc}
\usepackage{amssymb}
\usepackage{graphicx}
\begin{document}

\begin{center} {\bf How much information about a dynamical system do
its recurrences contain?}

{\bf Geoffrey Robinson$^1$, Marco Thiel$^2$}\\
{$^1$ Department of Mathematics, University of Aberdeen, AB24 3UE,
UK}\\
{$^2$ Department of Physics, University of Aberdeen, AB24 3UE, UK}\\
\end{center}

\begin{abstract}
We show that, under suitable assumptions, Poincar\'e recurrences of
a dynamical system determine its topology in phase space. Therefore,
dynamical systems with the same recurrences are topologically
equivalent.
\end{abstract}

\noindent{\bf Introduction:} From the beginning of what is today
called dynamical systems theory, there has been great interest in
the following problem: what can be inferred about a system if only
partial information about its state at each moment in time is given?
This matter is highly relevant in physics, where often not all state
variables of a system but only one \emph{observable} can be
measured. To make the following arguments more precise we introduce
the following definitions.

\medskip
\noindent {\bf Definition 1.1.} \emph{ Let $M$ be a compact
manifold. A dynamical system on M is a diffeomorphism
$\phi:M\longrightarrow M$ (discrete time) or a vector field $X$ on
$M$ (continuous time). In both cases the time evolution
corresponding to an initial position $x_0\in M$ is denoted by
$\phi_t(x_0)$. In the discrete case $t\in \mathbb{N}$ and
$\phi_i=(\phi)^i$. In the case of continuous $t\in\mathbb{R}$ and
$t\longrightarrow\phi_t(x_0)$ is the $X$ integral curve through
$x_0$.}\\

\medskip
\noindent {\bf Definition 1.2 :} An \emph{observable} is a smooth function $y:M\longrightarrow\mathbb{R}$.\\

\medskip
In the late 70's it was known that a procedure of delay embedding
could be used to reconstruct the \emph{"Geometry [of an attractor]
from a Time Series"} \cite{Packard79}. This reconstruction could
then be used to estimate dynamical invariants such as Lyapunov
exponents. In 1981, Takens proved three theorems \cite{Takens}, the
first of which reads (using definitions 1.1. and 1.2.)

\medskip
\noindent {\bf Theorem 1.1: } \emph{
Let $M$ be a compact manifold of dimension $m$. For pairs\\
$(\phi,y),\phi:M\longrightarrow M$ a smooth diffeomorphism and
$y:M\longrightarrow M$ a smooth function, it is a generic property
that the map $\Phi_{(\phi,y)}:M\longrightarrow\mathbb{R}^{2m+1}$,
defined by
\begin{equation}\label{delayemb}
\Phi_{\phi,y}(x)=(y(x),y(\phi(x)),\ldots,y(\phi^{2m}(x)))
\end{equation}
is an embedding; here "smooth" means at least $C^2$.}

\medskip
As mentioned above, the delay embedding procedure (\ref{delayemb})
lies at the heart of time series analysis and several authors have
extended Takens' Theorem \cite{Sauer,Stark}. These theorems are
virtually always used when nonlinear time series analysis is applied
to experimental data \cite{KantzSchreiber}. As the reconstructed
attractor is diffeomorphic to the original attractor, it can be
considered to be the same system after a smooth co-ordinate change.
Such a change leaves dynamical invariants, such as entropies,
fractal dimensions and Lyapunov exponents
unchanged.\\

In this paper we study a similar situation, in which it is not an
observable in the sense of definition 1.2. which is known, but
rather the Poincar\'e recurrences of the system. The theorems we
report are not primarily intended for applications in time series
analysis, but rather to stress the importance of Poincar\'e
recurrences. However, there are some applications, e.g. in the field
of data analysis namely for the method of recurrence plots
\cite{MarwanPhysRep} or for the reconstruction of protein structures
\cite{bohrprot} for which the results are relevant.

\medskip
\begin{center}
{\bf Poincar\'e recurrence, recurrence matrix and recurrence plots}
\end{center}

Recurrence is a fundamental characteristic of many dynamical
systems. The concept of recurrence was introduced by Poincar\'e in a
seminal work in 1890 \cite{poincare1890}, for which he won a prize
awarded by King Oscar~II of Sweden and Norway. Therein, Poincar\'e
not only described the ``homoclinic tangle'' which lies at the root
of the chaotic behaviour of orbits, but he also introduced (as a
by-product) the concept of recurrences in conservative systems. When
speaking about the restricted three body problem he reported :`` In
this case, neglecting some exceptional trajectories, the occurrence
of which is infinitely improbable, one can show, that the system
recurs infinitely many times as close as one wishes to its initial
state.'' In a larger context recurrences make part of one of three
broad classes of asymptotic invariants \cite{katok95}: (i) growth of
the number of orbits of various kinds and of the complexity of orbit
families \footnote{An important invariant of the orbit growth is the
topological entropy}, (ii) different types of recurrences, and (iii)
asymptotic distribution and statistical behaviour of orbits. The
first two classes are of purely topological
nature. The last class is naturally related to ergodic theory.\\
Of the different types of recurrences which form part of the second
class of invariants, the Poincar\'e recurrence is of particular
interest to this work. It is based on the Poincar\'e Recurrence
Theorem \cite{katok95} (Theorem 4.1.19):

\medskip
\noindent {\bf Theorem 2.1.} \emph{ Let $T$ be a measure-preserving
transformation of a probability space $(\mathcal{X},\mu)$ and let
$\mathcal{A} \subset \mathcal{X}$ be a measurable set\footnote{Here
$\mu$ is a Borel measure on a separable metrisable space
$\mathcal{X}$. Note that these assumptions are rather weak from a
practical point of view. Such a measure preserving function is
obviously given in Hamiltonian systems and also for all points on
(the $\omega$-limit set of) a chaotic attractor.}. Then for any
natural number $N \in \mathbb{N}$
\begin{equation}\label{eq_poincareTheorem}
\mu \Bigl(\bigl\{x\in \mathcal{A} \, |\  \{T^n(x)\}_{n\ge N}\subset
\mathcal{X} \backslash \mathcal{A}\bigr\} \Bigr)=0.
\end{equation}\vspace{1.0cm} }

Due to the great relevance of this theorem for this paper we give
here a rather short (standard) proof of this theorem:\\
\medskip
\noindent {\bf Proof:} Replacing $T$ by $T^N$ in
Eq.~(\ref{eq_poincareTheorem}), we find that it is enough to prove
the statement for $N=1$. The set
\[
\tilde{\mathcal{A}} := \left\{ x \in \mathcal{A} \, |\
\{T^n(x)\}_{n\in\mathbb{N}}
        \subset \mathcal{X} \backslash \mathcal{A}\right\} = \mathcal{A}
        \cap \left(\bigcap\limits_{n=1}^{\infty}T^{-n}(\mathcal{X} \backslash \mathcal{A})\right)
  \nonumber
\]
is measurable. Note that
\[
T^{-n}(\tilde{\mathcal{A}})\cap\tilde{\mathcal{A}}=\emptyset
\nonumber
\]
for every
$n$, because if we suppose, on the contrary, that $T^{-n}(\tilde{\mathcal{A}})\cap\tilde{\mathcal{A}}=B$ with
 $B \ne \emptyset$, this implies $T^{n}(B) \subset \mathcal{A}$. This is inconsistent with the definition
 of $\tilde {\mathcal{A}}$, because $B \subset \tilde{\mathcal{A}}$.\\
Also note that
\[
T^{-n}(\tilde{\mathcal{A}})\cap
T^{-m}(\tilde{\mathcal{A}})=\emptyset \quad \forall\,
m,n\in\mathbb{N} \nonumber,
\]
because if we assume, on the contrary, that
$T^{-n}(\tilde{\mathcal{A}})\cap T^{-m}(\tilde{\mathcal{A}})=B$ with
$B \ne \emptyset$, this implies $T^n(B) \subset \tilde{\mathcal{A}}$
and $T^m(B) \subset \tilde{\mathcal{A}}$. Without loss of
generality, we assume that $m>n$. Then
$T^n(B)=C \subset \tilde{\mathcal{A}}$ and $T^m(B)=T^{m-n}(T^n(B))=T^{m-n}(C)\subset \tilde{\mathcal{A}}$, which is again
inconsistent with the definition of $\tilde{\mathcal{A}}$.\\

Furthermore,
$\mu(T^{-n}(\tilde{\mathcal{A}}))=\mu(\tilde{\mathcal{A}})$ since
$T$ preserves $\mu$. Thus $\mu(\tilde{\mathcal{A}})=0$ because
\[
1 =
\mu(\mathcal{X})\ge\mu\left(\bigcup^{\infty}_{n=0}T^{-n}(\tilde{\mathcal{A}})\right)
=
   \sum^{\infty}_{n=0}\mu\left(T^{-n}(\tilde{\mathcal{A}})\right) =
   \sum^{\infty}_{n=0}\mu(\tilde{\mathcal{A}}).\qquad \Box\nonumber
\]
That means, that if we have a measure preserving transformation, the
trajectory will eventually come back to the neighbourhood of any
former point with probability one. However, the theorem only
guarantees the existence of recurrence but gives no indication of
how long it takes the system to recur. Especially for high
dimensional complex systems the recurrence time might be
extremely long.\\
The {\it first return} of a {\it set} is defined as follows: if
$\mathcal{A}\subset \mathcal{X}$ is a measurable set of a measurable
(probability) dynamical system $\left\{\mathcal{X},\mu,T\right\}$,
the first return of the set $\mathcal{A}$ is given by
\[
\tau (\mathcal{A})=\min\left\{n>0:T^n \mathcal{A} \cap \mathcal{A}
\not= \emptyset\right\}.
\]
Generically, for hyperbolic systems the recurrence or first return
time appears to exhibit certain universal properties
\cite{balakrishnan2001}:
\begin{enumerate}
\item the recurrence time has an exponential limit distribution; \label{consequ_1}
\item successive recurrence times are independently distributed; \label{consequ_2}
\item as a consequence of (\ref{consequ_1}) and (\ref{consequ_2}), the sequence of successive recurrence times has a
limit law that is a Poisson distribution.
\end{enumerate}
These properties, which are also well-known characteristics of
certain stochastic systems, such as finite aperiodic Markov chains
\cite{feller49,cox94,pitskel91}, have been rigorously established
for deterministic dynamical systems exhibiting sufficiently strong
mixing \cite{sinai70,hirata93,collet92}. They have also been shown
valid for a wider class of systems which have to be hyperbolic
\cite{hirata95}.\\
Recently, recurrences and return times have been studied with
respect to their statistics \cite{hirata99,penne99} and have been
linked to various other basic characteristics of dynamical
systems, such as the Pesin's dimension \cite{afraimovich97}, the
point-wise and local dimensions
\cite{afraimovich2000,afraimovich2003,gao99} or the Hausdorff
dimension \cite{barreira2001}. Also the multi-fractal properties
of return time statistics have been studied
\cite{hadyn2002,saussol2003a}. Furthermore, it has been shown that
recurrences are related to Lyapunov exponents and to various
entropies \cite{saussol2002,saussol2003b}. They have been linked
to rates of mixing \cite{young99}, and the relationship between
the return time statistics of continuous and discrete systems has
been investigated \cite{balakrishnan2000}.\\
An alternative and very practical approach has been introduced by
Eckmann et al. in 1987 \cite{Eckmann87}. Instead of studying the
(first) return time statistics, they introduced the recurrence
matrix
\begin{equation}\label{RPmatrix}
R_{i,j}=\Theta\left(\varepsilon-\left||\phi_i(x_0)-\phi_j(x_0)\right||\right),
\end{equation}
where $\phi_i(x_0),\phi_j(x_0)$ are the points on the continuous or
discrete orbit of a dynamical system at discrete times $i$ and $j$,
i.e. Eckmann et al. considered vector valued time series. If the
(Euclidean) distance between these two points is smaller than the
(arbitrary but small) threshold distance $\varepsilon$, the
Heaviside function $\Theta(\cdot)$ leads to $R_{i,j}=1$. If their
distance is larger than $\varepsilon$ we have $R_{i,j}=0$. Note that
we will adapt the definition of the Heaviside function for which
$\Theta(0)=0$. In principle, the definition $\Theta(0)=1$ would also
be possible. Based on these definitions, $R_{i,j}$ is a binary
matrix, which contains information about the recurrences
of the dynamical system.\\
This means that the orbit of the dynamical system is mapped to a
binary matrix, which indicates that, in a sense, much of the
information these orbits contained is lost. However, the recurrence
matrix and the derived recurrence plots, which are a graphical
representation of the recurrence matrix, have been applied very
successfully to analyse time series. From the matrix many dynamical
invariants, such as generalised Renyi entropies and dimensions and
also the mutual information, could be estimated \cite{ThielChaos}.
Recently, it has been shown that for one-dimensional systems the
time series can be reconstructed from the recurrence matrix up to a
smooth observation function \cite{Thiel_PLA}. Later a numerical
algorithm was used to reconstruct the attractor from the recurrence
matrix of various systems, e.g. the Lorenz attractor
\cite{Thiel_07}. An interesting feature of these algorithms is that
they are based on purely geometrical considerations, which is
different
from delay embedding procedure in Takens' Theorem. \\
These results suggested that once the recurrences of a system are
known (in the form of the recurrence matrix in the limit of long
time series) its dynamics are determined.\\
In the next section we will prove a theorem, which will show that
under rather weak assumptions it is possible to "reconstruct the
geometry" of an attractor if the {\em distance inequalities} are
known for each pair of points.

\section{Proof of the Reconstruction Theorem}

To investigate how much recurrences tell us about the dynamics of a
system, we suppose that we have a subset $M$ of $\mathbb{R}^{n}$ for
some $n,$ and we know for each pair of points their distance
inequality, i.e. we know if $d(x,y) < \varepsilon,$ or if $d(x,y)
\geq \varepsilon$. In terms of dynamical systems theory, these might
be the points on the attractor in phase space, together with their
"recurrence matrix". Suppose we then map all points of $M$ onto
another set $N$ by a function $f$, such that the resulting distance
inequalities are the same, i.e. we have $d(f(x),f(y)) < \varepsilon$
if and only if $d(x,y) < \varepsilon$. Now we have two sets (with an
induced topology inherited from $\mathbb{R}^{n}$), which have the
same distance inequalities, or from the point of view of dynamical
systems theory, the same recurrences. In the context of the
dynamical considerations, we impose a fairly natural additional
condition: that any two elements of $M$ can be distinguished by
comparing their distances (relative to $\varepsilon$) from enough
elements of $M.$ More precisely, if the ``separation'' condition of
the Theorem (and its Corollary) below fails for $M$, then there are
distinct points $x,y \in M$ such that the closed ball (in $M$) of
radius $\varepsilon$  with centre $x$ is contained in the
corresponding closed ball in $M$ with centre $y$. In the following
proofs, we need to distinguish carefully between strict inequalities
and cases where equality might be allowed. In the practical
applications, this distinction seems unlikely to matter.

A corollary, which is a consequence of the following theorem will
show, that $f$ must then
induce a homeomorphism.\\

\medskip
\noindent {\bf RECONSTRUCTION THEOREM:} \emph{ Let $M,N$ be closed
subsets of $\mathbb{R}^{n}$ for some $n < \infty.$ Let $f: M \to N$
be a surjection with the property that for some fixed $\varepsilon >
0,$ we have $d(f(x),f(y)) < \varepsilon$ if and only if $d(x,y) <
\varepsilon,$ where $d$ is the usual metric. Suppose also that, for
each $x\neq y \in M,$ there is an element $z \in M$ with $d(x,z)
> \varepsilon$ and $d(y,z) < \varepsilon,$ and similarly, whenever
$u \neq v \in N,$ there is an element $w \in N$ with $d(u,w)
> \varepsilon$ and $d(v,w) < \varepsilon.$ Then $f$ is a homeomorphism.}

\medskip
\noindent {\bf PROOF:} Notice that $f$ is injective, since for
$x,y,z$ above, $d(f(x),f(z)) \geq  \varepsilon$ and $d(f(y),f(z)) <
\varepsilon,$ so $f(x)$ and $f(y)$ are distinct. Hence, by symmetry,
it suffices to prove that $f$ is continuous.

\medskip
Notice, then, that for each $x \in M,$
$$M \backslash \{x\} = \bigcup_{\{z \in M: d(x,z) > \varepsilon \}}
B(z,\varepsilon),$$ where $B(z,\varepsilon)$ denotes the open ball
in $M$ with centre $z$ and radius $\varepsilon,$ since $x$ is
outside this union, but, by hypothesis, any $y$ different from $x$
is in the union (it is at distance less than $\varepsilon$ from some
$z$ such that $d(x,z) > \varepsilon ).$

Also,
$$N \backslash \{f(x)\} = \bigcup_{\{z \in M : d(x,z) > \varepsilon \}}
B'(f(z),\varepsilon),$$ where $B'(v,\varepsilon)$ denotes the open
ball in $N$ with  centre $v$ and radius $\varepsilon .$ To see this,
we apply $f$ to each side of the equality
$$M \backslash \{x\} = \bigcup_{\{z \in M: d(x,z) > \varepsilon \}}
B(z,\varepsilon)$$ and note that $$f(B(z,\varepsilon)) =
B'(f(z),\varepsilon)$$ by the hypotheses.

\medskip
Take a sequence $(x_{k})$ of elements of $M$ having limit $x \in M.$
We will prove that $(f(x_{k}))$ has limit $f(x).$  Let $K =
\overline{B'(f(x),\varepsilon)}$ (since $N$ is closed in
$\mathbb{R}^{n},$ the closures in $N$ and $\mathbb{R}^{n}$
coincide). Notice that $K$ is closed and bounded, so is compact. For
large enough $k,$ we have $x_{k} \in B(x,\varepsilon)$ and $f(x_{k})
\in B'(f(x),\varepsilon),$ using the hypotheses for the latter
inclusion.

Choose $\delta > 0.$ Now $$N  = B'(f(x),\delta) \cup \bigcup_{\{z
\in M : d(x,z) > \varepsilon \}} B'(f(z),\varepsilon),$$ and
$$K  = (K \cap B'(f(x),\delta)) \cup \bigcup_{\{z
\in M : d(x,z) > \varepsilon \}} (K \cap B'(f(z),\varepsilon)),$$
 a cover of $K$ by open sets (in $K$). Since $K$ is compact, there is a finite
subcover, which must include $(K \cap B'(f(x),\delta)),$ as none
of the open sets $B'(f(z),\varepsilon)$ (where $d(x,z) >
\varepsilon$) contains $ f(x),$ (for $d(x,z) > \varepsilon$ yields
$d(f(x),f(z)) \geq \varepsilon$ by hypothesis).

\medskip
Hence there is an integer  $s$ and elements $v_{1},\ldots, v_{s} \in
M$ such that $d(x,v_{j}) > \varepsilon$ and $d(f(x),f(v_{j}) \geq
\varepsilon$ (for each $j$) and
$$ K  = ( K \cap B'(f(x),\delta))  \cup
\bigcup_{j=1}^{s}(K \cap B'(f(v_{j}),\varepsilon)).$$ Let $$\alpha =
{\rm min} \{d(x,v_{j}) - \varepsilon : 1 \leq j \leq s\}
> 0.$$ For all large enough $m,$ we have $x_{m} \in B(x,\varepsilon),$
$f(x_{m}) \in B'(f(x),\varepsilon)$ and $d(x,x_{m}) < \alpha.$  Then
(for each $j$) $$d(x,v_{j}) \leq d(x,x_{m}) + d(x_{m},v_{j}) <
\alpha + d(x_{m},v_{j}),$$ so that $$d(x_{m},v_{j}) > d(x,v_{j})
-\alpha \geq d(x,v_{j}) -(d(x,v_{j}) - \varepsilon) = \varepsilon,$$
as $\alpha \leq d(x,v_{j}) - \varepsilon .$ Hence
$d(f(x_{m}),f(v_{j}) \geq \varepsilon$ for each $j.$

\medskip
Since $$ K = ( K \cap B'(f(x),\delta)) \cup \bigcup_{j=1}^{s}( K
\cap B'(f(v_{j}),\varepsilon))$$ and $d(f(x_{m}),f(v_{j})) \geq
\varepsilon$ for each $j,$ it must be the case that
$d(f(x),f(x_{m})) < \delta.$ Since $\delta$ is arbitrarily chosen,
$(f(x_{k}))$ has limit $f(x).$ Thus $f$ is continuous, as required.
$\Box$\\

The following useful corollary is immediate, since compact subsets
of $\mathbb{R}^{n}$ are closed, and because the separation property
required in $N$ follows automatically from the properties
of $f$ and the fact that the corresponding property holds in $M.$\\

\medskip
\noindent {\bf Corollary 1:} \emph{ Let $M,N$ be compact subsets of
$\mathbb{R}^n$ for some $n < \infty$. Let $f : M \longrightarrow N$
be a surjection with the property that for some fixed $\varepsilon >
0$, we have $d(f(x), f(y)) < \varepsilon$ if and only if $d(x, y) <
\varepsilon$ and $d(f(x), f(y))
> \varepsilon$ if and only if $d(x, y) > \varepsilon$, where $d$ is the usual metric. Suppose
also that, for each $x \ne y \in M$, there is an element $z \in M$
with $d(x, z) > \varepsilon $ and $d(y, z) < \varepsilon$. Then f
is a homeomorphism.}

\begin{center}
{\bf Interpretation}
\end{center}

The theorem and the corollary given in the last section show that
distance inequalities restrict the possible structure of a set up to
homeomorphic transformations. The recurrence matrix also contains
information on the time evolution of the system. The indices $i$ and
$j$ in \Eref{RPmatrix} correspond to the times when the points
$\phi_i(x_0)$ and $\phi_j(x_0)$ are close to each other. This means
that the recurrence matrix contains the information of whether or
not the system at time $j$ has recurred to the state it was in at
$i$. If the system lives in its phase space $M$, the theorem and the
corollaries guarantee that if we take all the recurrences into
account, we can reconstruct a homeomorphic set $N$, once points with
identical neighbourhoods are identified. We furthermore know the
sequence in time in which the different points in the phase spaces
are visited. We therefore can reconstruct a topologically conjugate
dynamical system from the knowledge of the recurrences. This means
that the reconstructed system is dynamically conjugate to the
original one, in the sense of \cite{Robinson}. Both systems can then
be considered to be identical up to a continuous change of the
coordinate system.

\begin{center}
{\bf Conclusions}
\end{center}

Several conclusions can be drawn from these considerations. First of
all, the theorem shows that under some assumptions recurrences
determine a system up to topological conjugacy. Two systems with the
same recurrence matrix are dynamically equivalent and can be
considered to be linked by a continuous change of coordinates. This
result also sheds light on standard techniques which are used to
detect generalised synchronisation, where one basically considers
simultaneous recurrences in both systems \cite{GS1,Romano}. These
methods quantify the similarity of both recurrence matrices for
small thresholds $\varepsilon$, yielding an index for generalised
synchronisation. Furthermore, the main result of this paper provides
a theoretical basis for the numerical algorithm used in
\cite{Thiel_07} to reconstruct the trajectory of the underlying
system from its recurrence matrix. This reconstruction was shown to
be possible for a rather large interval of values of the threshold
$\varepsilon$. In the case of the Bernoulli map, for example, the
$\varepsilon$-ball could cover up to half of the phase space. Also
for continuous systems like the R\"ossler or the Lorenz system this
reconstruction could be performed for $\varepsilon$-neighbourhoods
covering a large fraction of the phase space. The quality of the
reconstruction was practically independent of $\varepsilon$ then, as
even spatial structures much smaller than $\varepsilon$
could be reconstructed for long time series.\\
The theorem proved in this paper also supports the method proposed
in  \cite{ThielTwin} to generate twin surrogate data from the
recurrence matrix.\\
Furthermore, our considerations suggest that it be might possible to
show that $f$ together with its inverse $f^{-1}$ has to be
differentiable, i.e. at least $\mathcal{C}^1$, and therefore induces
a diffeomorphism between $M$ and $N$.\\
This work illustrates that, in a sense, recurrences provide an
alternative description of a dynamical system.

\noindent {\bf Acknowledgement:}  The authors want to thank  M.
Carmen Romano, Celso Grebogi, J\"urgen Kurths,
Jaroslav Stark and Valentin Afraimovich for helpful and stimulating discussions.\\

\end{document}